\documentclass%[12pt]
    {article}
\usepackage{theorem}
\usepackage{oldgerm}
\usepackage{xspace}
\usepackage{enumerate}
\usepackage{latexsym}
\usepackage{amssymb}
\usepackage{amsmath}
\usepackage{epsfig}

\theoremstyle{plain}
        \newtheorem{thm}{Theorem}[section]
        \newtheorem{lem}[thm]{Lemma}
        \newtheorem{conj}[thm]{Conjecture}
        
        \newtheorem{prop}[thm]{Proposition}
        \newtheorem{cor}[thm]{Corollary}
        
        \newtheorem{probl}[thm]{Problem}
        {\theorembodyfont{\rmfamily}
                \newtheorem{defn}[thm]{Definition}
                \newtheorem{rem}[thm]{Remark}
                \newtheorem{exa}[thm]{Example}
        }

\def\negyzet{\vbox{\hrule
                   \hbox{\vrule\kern2pt
                         \vbox{\kern2pt\kern2pt
                        }\kern2pt\vrule
                  }\hrule
                   }
             }

\renewcommand{\rho}{\varrho}
\renewcommand{\phi}{\varphi}

\newcommand{\proof}{\par\noindent {\sl Proof: \/}}

\newcommand{\proofend}{~\rule{2mm}{3mm}}
\newcommand{\rnh}{$\le\rho$--neighborhood\xspace}

\newcommand{\utacfC}{up to a constant factor $C$\xspace}

\newcommand{\EP}{Erd\H{o}s\xspace}

\newcommand{\defeq}{\buildrel{\scriptstyle\rm def}\over=}

\def\dist{\hbox{\text{dist}}}

\newcommand{\ds}{\textstyle}

\newcommand{\wfp}{with the following property\xspace}

\def\CP{Cartesian product}
\def\SzTr{Szemer\'edi--Trotter\xspace}

\def\pni{\par\noindent}
\newcommand{\p}{^{\prime}}
\newcommand{\pp}{^{\prime\prime}}
\newcommand{\ppp}{^{\prime\prime\prime}}
\newcommand{\pppp}{^{\prime\prime\prime\prime}}

\newcommand{\A}{{\mathcal A}}

\newcommand{\Cell}{{\mathcal C}}

\renewcommand{\P}{{\mathcal P}}

\newcommand{\R}{\mathbb R}

\renewcommand{\S}{\mathcal S}

\newcommand{\U}{{\mathcal U}}

\begin{document}

\date{}

\title {On the Dimension of Finite Point Sets II.\\
        ``Das Budapester Programm''}
\author{Gy\"orgy Elekes}
\maketitle

\begin{abstract}\noindent

\end{abstract}

\begin{section}{Introduction}

\begin{subsection}{History: ``Das Erlanger Programm'' of Felix Klein}

``Klein's synthesis of geometry as the study of the properties 
of a space that are invariant under a given group of transformations, 
known as the Erlanger Programme (1872), profoundly influenced 
mathematical development. The Erlanger Programme gave a 
unified approach to geometry which is now the standard 
accepted view.''
(Article by: J.~J.~O'Connor and E.~F.~Robertson)
\pni
{\tt
%http://
www-history.mcs.st-andrews.ac.uk/history/Mathematicians/Klein.html}
\end{subsection}

\begin{subsection}{A modest sub-program.}

We are going to use the nickname ``Das Budapester Programm''
for a large class of finite combinatorial problems
to be posed below.
Following Klein's original idea, they concern
 various dimensional Euclidean spaces and  several
groups of transformations.
\par
Of course, we do not expect unifying any theory,
we just  hope that the questions may be interesting on their own right.
The one dimensional special cases (which are already solved) 
have non-trivial applications both to algebraic and geometric problems
(see Section~\ref{OneDimSect}).
That is why we believe that also the higher dimensional versions
--- of which we solve one here --- may be useful.
\medskip\par
Given a group $G$ of transformations of $\R^d$ %(or $\C^d$)
and a finite point set $\P\subset\R^d$ %(or $\P\subset\C^d$),
we shall be interested in the number of transformations 
$\phi\in G$
which map many points of $\P$ to some other points of $\P$:
$$
\#\{\phi\in G\;;\; \phi(\P)\cap\P \text{ is ``not too small''} \}.
$$
Usually, ``not too small'' will mean $|\phi(\P)\cap\P|\ge k$
for a given (large) $k\le|\P|$.
However, sometimes we shall even require that $\phi(\P)\cap\P$
be ``proper $d$--dimensional'' in some sense,
 introduced in Part I \cite{EGy:DoFI}.
\begin{defn}
\label{ProperDimDef}
A set of $N$ points is \emph{proper $d$--dimensional 
\utacfC,\/} (for short, ``proper $d$--D'') 
if it  can be  cut into singletons by
 at most $ C\root{d}\of{N}$ appropriate hyperplanes.
\end{defn}
The ultimate goal of our ``Budapester Programm'' 
is to find sharp upper bounds (as a function of $N$,
for various dimensions $d$ and groups $G$) for
$$
f_G^d(N,k)\defeq
\max_{{{|\P|=N \atop\P\subset\R^d } \atop \text{proper }d\text{--dim}
       }}
\#\{\phi\in G\;;\; 
    \phi(\P)\cap\P  \text{ proper $d$--dim, } |\phi(\P)\cap\P|\ge k\}.
$$
\begin{probl}[``Meta--problem'']
Is it true that, for any ``reasonable Euclidean'' group $G$
which acts on $\R^d$, 
$$
f_G^d(N,k)=O\biggl(\frac{N^\alpha}{k^{\alpha-1}}\biggr),
$$
for an $\alpha=\alpha(G)>1$, independent from $N$ and $k$?
\end{probl}
\end{subsection}

\begin{subsection}{More history: old results in $\R$.}
\label{OneDimSect}

In one dimension, the assumption of being proper 1D
only means that the points of $\P$ are all distinct.
\begin{prop}\label{AffineProp}
For the group of non-constant affine transforms of $\R$,
i.e., for  $G=\{\phi:\R\rightarrow\R\;;\; \phi(x)=mx+b,\ m\ne0\}$,
we have
$$
\#\{\phi:\R\rightarrow\R \text{ affine }\;;\; |\phi(\P)\cap\P|\ge k\}
    \ \le\   C\cdot|\P|^4/k^3,
$$
for an absolute constant $C$ and any\/ $2\le k \le |\P|$.
\end{prop}
\proof
See \cite{EGy:LinI} and Section~\ref{OneDimFinePts} for some remarks;
also \cite{EGy:survey} for a survey.
\begin{prop}\label{ProjAndRatProp}
For the group of non-constant projective transforms of $\R$,
i.e., for  $G=\{\phi:\R\rightarrow\R\;;\; \phi(x)=(ax+b)/(cx+d),\ 
ac-bd\ne0\}$,
we have
$$
\#\{\phi:\R\rightarrow\R \text{ projective }\;;\; |\phi(\P)\cap\P|\ge k\}
    \ \le\   C\cdot|\P|^6/k^5,
$$
for an(other) absolute constant $C$ and any\/ $3\le k \le |\P|$.
\par
For any positive integer $r$ and for the group $G_r$ of 
non-constant rational transforms of $\R$ of total degree $\le r$,
i.e., for  $G=\{\phi:\R\rightarrow\R\;;\; \phi(x)=p(x)/q(x),
\ p,q\in\R[x],\; \deg(p)+\deg(q)\le r\}$,
we have
$$
\#\{\phi\in G_r\;;\; |\phi(\P)\cap\P|\ge k\}
    \ \le\   C\cdot|\P|^{2r+2}/k^{2r+1},
$$
for a constant $C=C(r,c)$.
\end{prop}
\proof
See \cite{EGyKZ}, \cite{EGyRL} and 
Section~\ref{OneDimFinePts} for some remarks.
\smallskip\par
The upper bound in  Proposition~\ref{AffineProp} gives
the best possible order of magnitude.
However, this is unknown for those in Proposition~\ref{ProjAndRatProp};
perhaps they can be improved for proper (i.e., non-affine)
rational or even projective transforms.
\begin{rem}
\label{LinBdRem}
It is worth noting that, for any fixed $0<c<1$ and $k=c|\P|$,
the bounds in Propositions~\ref{AffineProp} and \ref{ProjAndRatProp}
are linear in $|\P|$, i.e.,
$\#\{\phi:\R\rightarrow\R \text{ affine or projective or rational}\;;\; 
      |\phi(\P)\cap\P|\ge c|\P|\} \ \le\   C\cdot|\P|.
$
Though this only is a very special case, bounds like this 
were crucial in finding Freiman--Ruzsa type structure results for
small composition sets of affine or projective transforms of $\R$
in \cite{EGY:LinII,EGyKZ}.
\end{rem}
\end{subsection}
\end{section}

\begin{section}{Some problems and an affine result in $\R^2$.}

\begin{subsection}{Isometries.}
The first interesting case is the group of isometries of 
the Euclidean plane, i.e.,
$G=\{\phi:\R^2\rightarrow\R^2\;;\; 
            \phi \text{ preserves distances}\}$.
Even the following is unknown.
\begin{conj}
There is an absolute constant $C$ such that
$$
\#\{\phi:\R^2\rightarrow\R^2 \text{ isometry}\;;\; 
      |\phi(\P)\cap\P|\ge 2\} \ \le \ C\cdot|\P|^3.
$$
\end{conj}
This, if true, could provide a missing link to
the ``distinct distances'' problem of \EP.

\begin{rem}
If we restrict ourselves to translations (=shifts),
then it is not difficult to show --- even in arbitrary 
dimension --- that
$$
\#\{\phi:\R^d\rightarrow\R^d \text{ translation}\;;\; 
      |\phi(\P)\cap\P|\ge k\} \ \le \ C\cdot|\P|^2/k,
$$
for an absolute constant $C$ and any $1\le k\le|\P|$.
(Just draw the --- at least $k$ --- ``arrows''
which indicate which point is mapped to which other point.)
Moreover, this order of magnitude is best possible,
as shown by Example~\ref{ShiftExa}.
\end{rem}
\end{subsection}

\begin{subsection}{The affine group of $\R^2$.}

Our principal result Theorem~\ref{MainThm}
concerns the affine group of $\R^2$,
i.e. that of mappings
$\phi(x,y)=(a_1x+b_1y+c_1,\ a_2x+b_2y+c_2)$ where
$a_1b_2-a_2b_1\ne0$.
Equivalently --- using projective coordinates ---
 this is the group of non-singular matrices
\[
\begin{pmatrix}
a_1&b_1&c_1\cr
a_2&b_2&c_2\cr
0&0&1\cr
\end{pmatrix}
\]
We shall restrict ourselves to the case $|\phi(\P)\cap\P|\ge k=c|\P|$.
According to Remark~\ref{LinBdRem},
we may expect a bound linear in $|\P|$,
even here in the two dimensional plane.

However, we must face several degenerate situations.
First, if $\P$ is collinear then there are infinitely many
affine mappings $\phi$ with  $\phi(\P)=\P$ --- whence
 $|\phi(\P)\cap\P|=|\P|$.
Also requiring that  $\phi(\P)\cap\P$ be non-collinear
--- and, of course, $|\phi(\P)\cap\P|\ge c|\P|$ ---
will not help: there still are point sets with a
super--linear (actually, quadratic) number of such mappings
(see Example~\ref{NonCollinAffineExa}).
\par
That is why in our main result we must assume that 
$\phi(\P)\cap\P$ is proper 2--dimensional.
(Note that a positive proportion, i.e., $cN$ points of a 
proper 2--dimensional set 
is still proper 2--dimensional 
--- up to another (larger) constant factor $C/\sqrt{c}$.)
\begin{thm}[Main Theorem]
\label{MainThm}
Let $\P\subset\R^2$ be proper 2--dimensional 
\utacfC.
Moreover, let $0<c<1$ be arbitrary.
Then
$$
\#\{\phi:\R^2\rightarrow\R^2 \text{ affine }\;;\; 
      |\phi(\P)\cap\P|\ge c|\P|\} \ \le\   C^*\cdot|\P|,
$$
where $C^*=C^*(C,c)$ does not depend on $|\P|$.
\end{thm}
The following questions remain open.
\begin{conj}
For any $\P\subset\R^2$, any $3\le k\le|\P|$ and any $C>0$,
\[
\begin{aligned}
\#&\{\phi:\R^2\rightarrow\R^2 \text{ affine }\;;\; 
    \phi(\P)\cap\P \text{ proper 2D up to $C$, and }
      |\phi(\P)\cap\P|\ge k\}  \cr
&\le\   C^*\cdot|\P|^6/k^5.\cr
\end{aligned}
\]
\end{conj}
If true, this order of magnitude (as a function of $|\P|$ and $k$)
is best possible, as shown by Example~\ref{ManyProperAffineExa}.
\begin{probl}
Does Theorem~\ref{MainThm} hold if, instead of $\P$
being proper 2D, we only require that
$$
\frac{|\text{maximal collinear subset of }\P|}{|\P|}
  \rightarrow0 \text{\quad while } {|\P|}\rightarrow\infty?
$$
\end{probl}
\end{subsection}
\end{section}

\begin{section}{Higher dimensions}

\begin{conj}
\label{HighDimAffineConj}
In $\R^d$, for  the group $G$ of non-degenerate
(bijective) affine mappings,
$$
f_G^d(N,k)=O\biggl(\frac{N^{2d+2}}{k^{2d+1}}\biggr).
$$
\end{conj}
The case $d=2$ was our Main Theorem~\ref{MainThm}.
We mention without detailed proof that it can  also
be extended to $d=3$ (see Remark~\ref{ThreeDimProofRem}).
However, we know nothing about $d\ge4$.
%
%The order of magnitude conjectured, if true, is best possible 
%(see Remark~\ref{ManyProperAffineRem}).
%
\medskip\par
We believe that the main feature of our Main Theorem~\ref{MainThm}
is that we could prove a sharp bound  on the number of certain mappings
\emph{\/even without an incidence bound at hand.}
Of course, life would be much easier if we had such bounds, e.g.
the following.
\begin{conj}
\label{HighDimAffSubSpaceConj}
Let $C>0$ and $\P\subset\R^D$ proper $D$--dimensional up to $C$.
We consider affine subspaces $S^{(r)}\subset\R^D$ of dimension $r<D$.
$$
\#\{ S^{(r)}\subset\R^D  \;;\; S^{(r)}\cap\P \;
\text{is $r$--dim up to $C$ and } 
|S^{(r)}\cap\P|\ge k \} \;  \le \; C^*\frac{|\P|^{r+1}}{k^{D+1}}.
$$
\end{conj}
(Part I, Problem~4.1 is the special case $r=1$.)\\
This bound can, again, be attained, see 
Example~\ref{HighDimAffSubSpaceExa}.
\begin{rem}
This conjecture would immediately imply the previous one.
(Let $\P_0\subset\R^d$ be an $N$--element set which attains
$f_G^d(N,k)$ and apply Conjecture~\ref{HighDimAffSubSpaceConj}
with $D=2d$ and $r=d$
to the $N^2$--element set $\P\defeq\P_0\times\P_0$ 
and the graphs --- as subsets of $\R^d\times\R^d$ --- of the non-degenerate
affine mappings $\phi:\P_0\rightarrow\P_0$.)
\end{rem}
\end{section}

%\vfill\eject

\begin{section}{Fine points and point sets}

\begin{subsection}{Proof of Propositions~\ref{AffineProp} 
     and \ref{ProjAndRatProp}}
\label{OneDimFinePts}

As for Proposition~\ref{AffineProp}, 
it is enough to consider the Cartesian product
$\P\times\P\subset\R^2$ of $n=|\P|\times|\P|$ points and 
use an incidence bound of \SzTr \cite{SzTr:83} which states that,
for any $n$ points of $\R^2$ and any $2\le k\le\sqrt{n}$,
at most $C\cdot n^2/k^3$ straight lines can pass through
$\ge k$ points each.

The proof of Proposition~\ref{ProjAndRatProp}
is based upon the same idea, using (a special case of) 
a result of Pach and Sharir \cite{PaSh:98} which states that,
for any $n$ points of $\R^2$ and any $3\le k\le\sqrt{n}$,
at most $C\cdot n^s/k^{2s-1}$  members of a family of curves 
of $s$ degrees of freedom can pass through $\ge k$ points each.\\
On the one hand, for  hyperbolas $y=(ax+b)/(cx+d)$ 
--- which form a family of three degrees of freedom ---
this gives a bound of 
$C\cdot n^3/k^5 =C\cdot |\P|^6/k^5$.
On the other hand, the degree of freedom for rational functions
of total degree $r$ is $s=r+1$, resulting in an upper estimate
$C\cdot n^{r+1}/k^{2r+1} =C\cdot|\P|^{2r+2}/k^{2r+1}$.
\end{subsection}

\begin{subsection}{Examples in $\R^2$}

\begin{exa}\label{ShiftExa}
(Our original argument was probabilistic;
here we present a simplified version of a construction by G\'eza T\'oth.)
\par
First we select a set $Y=\{y_1,\ldots,y_t\}$ of  $t:=n/(2k)$ reals
--- considered as points on the $y$-axis ---
such that the differences $y_i-y_j$ are all distinct for $i\ne j$. 
Then we define our point set as the $2k\times t$ Cartesian product
$$
\P:=\{1,2,\ldots,2k\}\times Y\subset \R^2.
$$
Now, for each pair $i\ne j$, we have at least $k$ distinct 
translations which map $k$ or more points at level $y_i$ to points at level
$y_j$, yielding a total of
$$
t(t-1)\cdot k \approx 
\biggl(\frac{n}{2k}\biggr)^2\cdot k =
\frac{1}{4}\cdot\frac{n^2}{k}.
$$
% If $k\le cN$ for a $0<c<1$,
% select a set $\P\subset\R^2$ of $\;\approx N$ random points 
% from a $(cN/\sqrt{k}) \times (cN/\sqrt{k})$ square lattice $\L$,
% each lattice point selected with probability $p=k/(c^2N)$.
% This way we really expect $(c^2N^2/k)\cdot p = N$ points.
% With high probability, at least $c^*N^2/k$ translations $\phi$
% will satisfy $|\phi(\P)\cap\P|\ge k$.\\
% To see this, we consider the
% set $\T$ of at least $(c^2N^2)/(4k)$ translations
% which map the bottom left quarter of $\L$ into $\L$.
% The probability that both a point and its image is preserved is
% $p^2=k^2/(c^4N^2)$ thus we expect a $T\in\T$ to map
% $p^2\cdot(c^2N^2)/(4k)=k/(4c^2)$ points of $\P$ into $\P$. 
% This quantity is much larger than $k$ if $c$ is small enough,
% hence the probability that $|T(\P)\cap\P|<k$ is exponentially small
% as a function of $c$ (or, rather, of $1/c$).
\end{exa}

\begin{exa}\label{NonCollinAffineExa}
Let $\P$ consist of $N=2n$ points: $n$ of them on the $x$-axis
(but otherwise arbitrary) plus $n$ other points $(0,y_i)$ ($1\le i\le n$)
such that no $y_i$ is~0 and the quotients $y_i/y_j$ ($i\ne j$)
are all distinct.\\
Then each of the affine mappings $(x,y)\mapsto(x,\lambda y)$
($\lambda=y_i/y_j$; $i\ne j$) map all $n$ points on the $x$--axis
to themselves plus exactly one point on the $y$--axis to another 
such point --- a total of $n+1>N/2$. 
Moreover, the number of such mappings is
$n(n-1)=N(N-2)/4$, indeed quadratic in $N=|\P|$.
\end{exa}
%\begin{exa}
%Let $\P$ be the following set of $N=4n+2$ points:
%$$
%\P\defeq\{(i,0)\;;\;i=-n,-(n-1),\ldots,n-1,n\} \cup 
%    \{(j,1)\;;\;j=-n,-(n-1),\ldots,n-1,n\},
%$$
%\[
%\begin{pmatrix}
%1&j-i&i\cr
%0&1&0\cr
%0&0&1\cr
%\end{pmatrix}
%\]
%\end{exa}

\begin{exa}\label{ManyProperAffineExa}
Let $k\le n/4$ and put $t=\sqrt{k}$. 
Define the $t\times(n/t)$ \CP
$\;\P\defeq\{1\ldots t\}\times\{1\ldots n/t\}\subset\R^2$, 
a set  of $n$ elements.
Then there are $\sim n^6/k^5$ non-degenerate affine mappings $\phi$ 
for which $\phi(\P)\cap\P$ contains a parallelogram--lattice of
$k=t\times t$ points.
\end{exa}
\proof
$\P\times\P\subset\R^4$ is a $t\times(n/t)\times t\times(n/t)$ \CP.
Using $(x,y,z,w)$ as coordinates, we describe some planes in $\R^4$,
parameterized by $x$ and $z$ (NOT $y\;$!),
and we shall make sure that these planes be graphs of non-degenerate
(bijective) affine mappings $\phi:\R^2\rightarrow\R^2$.
\par
Let the planes be determined by the equations
\[
\begin{aligned}
y&=a_1 x+b_1 z +c_1 \cr
w&=a_2 x+b_2 z +c_2, \cr
\end{aligned}
\]
such that $a_1b_2-a_2b_1\ne0$ and
\[
\begin{aligned}
a_1,a_2,b_1,b_2\in&\biggl\{1,2,\ldots,\frac{1}{3}\cdot\frac{n}{t^2}\biggr\};\cr
c_1,c_2\in&\biggl\{1,2,\ldots,\frac{1}{3}\cdot\frac{n}{t}\biggr\}.\cr
\end{aligned}
\]
Then we have 
$$
\biggl(\frac{1}{3}\cdot\frac{n}{t^2}\biggr)^4\cdot
  \biggl(\frac{1}{3}\cdot\frac{n}{t}\biggr)^2=
   \frac{1}{729}\cdot\frac{n^6}{t^{10}}
$$
such planes, of which, given any triple $(a_1,b_1,a_2)$,
at most one $b_2$ can violate  $a_1b_2-a_2b_1\ne0$.
Hence, at least
$$
 \frac{1}{729}\cdot\frac{n^6}{t^{10}}\cdot
    \biggl(1-\frac{3t^2}{n}\biggr)\ge
    \frac{1}{729}\cdot\frac{n^6}{t^{10}}\cdot\frac{1}{4}>
   \frac{1}{3000}\cdot\frac{n^6}{t^{10}}
$$
planes satisfy the extra condition.
Moreover, each of these contains a $k=t\times t$ parallelogram lattice, 
e.g., the points which correspond to the values $x,z\in\{1,\dots,t\}$,
since then
$$
y,w=a_ix+b_iy+c_i\le
   \biggl(\frac{1}{3}\cdot\frac{n}{t^2}\biggr)\cdot t +
   \biggl(\frac{1}{3}\cdot\frac{n}{t^2}\biggr)\cdot t +
   \frac{1}{3}\cdot\frac{n}{t} = \frac{n}{t}.
$$
We are left to show that these planes are graphs of bijective 
affine mappings.
This follows from re-writing the equations as
\[
\begin{aligned}
z=-\frac{a_1}{b_1}&\cdot x+\frac{1}{b_1}\cdot y-\frac{c_1}{b_1}\cr
w=\frac{a_2b_1-a_1b_2}{b_1}&\cdot x+\frac{b_2}{b_1}\cdot y+
                                 \frac{c_2b_1-c_1b_2}{b_1}\cr
\end{aligned}
\]
whose determinant is
$$
\det=-\frac{a_1}{b_1}\cdot\frac{b_2}{b_1} -
  \frac{1}{b_1}\cdot\frac{a_2b_1-a_1b_2}{b_1}=\frac{a_2}{b_1}\ne0,
$$
since $a_2\ne0$.
We conclude that we have  $(1/3000) n^6/k^5$ non-degenerate affine mappings 
$\phi$ for which $\phi(\P)\cap\P$ is proper 2--dimensional, as required,
since it contains a parallelogram--lattice of $k=t\times t$ points.
\proofend
\end{subsection}

\begin{subsection}{An example in $\R^D$}
%\begin{rem}
%\label{ManyProperAffineRem}
%Similar construction, at least for $d$ even ?!?????!?!?!
%\end{rem}
\begin{exa}
\label{HighDimAffSubSpaceExa}
Let $r<D$ and $k\le N^{r/D}$ be given and assume 
(at the cost of a constant factor)  
that $k=t^r$ for a positive integer $t\le N^{1/D}$.
Define $\P\subset\R^D$ as a cube lattice of size
$$
t\times t\times \ldots \times t \times 
   \biggl(\frac{N}{k}\biggr)^{\frac{\ds 1}{\strut\ds D-r}}
         \times\ldots\times
   \biggl(\frac{N}{k}\biggr)^{\frac{\ds 1}{\strut\ds D-r}}
$$
a product of $r+(D-r)$ terms.
Consider the $r$--dimensional affine subspaces determined by the
systems of equations
\[
\begin{aligned}
x_{r+1} =& a_1^{(r+1)} x_1 + \ldots + a_r^{(r+1)} x_r + c^{(r+1)}   \cr
%.&\cr
\vdots&\cr
%.&\cr
x_{D} =& a_1^{(D)} x_1 + \ldots + a_r^{(D)} x_r + c^{(D)},   \cr
\end{aligned}
\]
for
$$
a_i^{(.)}=1,2,\ldots,
  \frac{1}{t(r+1)}\cdot\biggl(\frac{N}{k}\biggr)^{\frac{\ds 1}{\strut\ds D-r}}
\text{ \quad and \ }
c^{(.)}=1,2,\ldots,
  \frac{1}{r+1}\cdot\biggl(\frac{N}{k}\biggr)^{\frac{\ds 1}{\strut\ds D-r}}.
$$
Then each such affine subspace contains $k=t\times t\times \ldots\times t$
points of $\P$, i.e. those for $x_1,x_2,\ldots,x_r\in\{1,2,\ldots,t\}$.
Moreover, the number of them is at least
\[
\begin{aligned}
\ge&\Biggl[
\frac{1}{t(r+1)}\cdot\biggl(\frac{N}{k}\biggr)^{\frac{\ds 1}{\strut\ds D-r}}
\Biggr]^{(D-r)\cdot r}
\cdot
\Biggl[
  \frac{1}{r+1}\cdot\biggl(\frac{N}{k}\biggr)^{\frac{\ds 1}{\strut\ds D-r}}
\Biggr]^{(D-r)} \ \approx \cr
&\ \approx c(r,D)\cdot 
      \biggl(\frac{1}{k}\biggr)^{D-r} 
	\cdot\biggl(\frac{N}{k}\biggr)^{r+1} = 
    c(r,D)\cdot\frac{N^{r+1}}{k^{D+1}},\cr
\end{aligned}
\]
where $c(r,D)>0$ does not depend on $N$ or $k$.
\end{exa}
\end{subsection}
\end{section}

%\vfill\eject

\begin{section}{Some Lemmata and the main proof}

\begin{subsection}{Arrangements of straight lines in $\R^2$.}

Let $H$ be a (finite) set of $n$ planes in $\R^2$.
They cut the plane into at most 
${n \choose 2} + {n \choose 1} + {n \choose 0} \sim n^2$
open convex cells,
with equality iff $H$ is in general position,
i.e., if any two intersect but no three do.

The set of cells, together with their vertices and edges,  
is called the \emph{arrangement\/} defined by $H$. 
We shall denote it by $\A(H)$.

For two cells $\Cell_i,\Cell_j\in \A(H)$, a natural notion 
of distance is
$$
\dist(\Cell_i,\Cell_j)\defeq\#\{h\in H\ ;\ h \text{ separates }	
		\Cell_i \text{ and } \Cell_j\}.
$$
It is easy to see that ``dist'' is a metric, 
i.e. it satisfies the triangle inequality.

The result presented below bounds --- from above, in terms of $|H|$ ---
the number of pairs $(\Cell_i,\Cell_j)$ whose distance is
at most a given $\rho>0$.
To this end, we also defined in Part I \cite{EGy:DoFI}
 the \rnh of a cell $\Cell_j$ as
$$
B_{\rho}(\Cell_j) = \{\Cell_i\in \A(H)\ ;\ 
		\dist(\Cell_i,\Cell_j)\le \rho	\},
$$
and the number of ``$\rho$--close pairs'' mentioned above equals 
$\sum
   %_{\Cell_j\in\A(H)}	
     |B_{\rho}(\Cell_j)|$.
\begin{prop}
\label{PlanarBdProp} 
%We have 
$\sum_{\Cell_j\in\A(H)}	|B_{\rho}(\Cell_j)|= O(\rho^2|H|^2)$ 
in $\R^2$ and 
%in $\R^2$ while
%$\sum_{\Cell_j\in\A(H)}|B_{\rho}(\Cell_j)|= $
$=O(\rho^3|H|^3)$ in $\R^3$.
\end{prop}
\proof
See Part I \cite{EGy:DoFI} , Corollary 3.3 (??) and Lemma 3.7 (??).
\medskip\par
As for the second moment $\sum|B_{\rho}(\Cell_j)|^2$, 
it may not always be bounded by a quadratic function of $|H|$
(e.g., if the lines all surround a regular polygon).
\begin{probl}
\label{RefinementProbl}
Let $\A(H)$ be a simple  arrangement in $\R^2$.
Is it true that it can be refined to an $\A(H^+)$
by adding $O({|H|})$ new straight lines such that
$\sum_{\Cell_j\in\A(H^+)}|B_{\rho}(\Cell_j)|^2= O(\rho^4|H|^2)$?
\end{probl}
It may well be true that one can even force the stronger upper bound
$|B_{\rho}(\Cell_j)|=O(\rho^2)$ for each $\Cell_j\in\A(H^+)$
--- but it is ``even more'' unknown.
However, the following weaker version is an easy consequence of 
Proposition~\ref{PlanarBdProp}.
\begin{cor}
\label{UniformBallsCor}
There is an absolute constant $K_1$ \wfp.\\
Let $\A(H)$ be an arrangement, $\rho\le|H|$ and $c>0$ arbitrary.
Given any subset $A_0\subset\A(H)$ of \ $\ge c|H|^2$ cells,
there is a sub-subset $A_1\subset A_0$ of at least $(c/2)|H|^2$ cells
(i.e., half of them) which satisfy
$$
|B_{\rho}(\Cell_j)|\le\frac{K_1}{c}\rho^2 \text{ for each }\Cell_j\in A_1.
$$
\end{cor}
\proof
Let $K_0$ be the absolute constant hidden in the right hand side
of Proposition~\ref{PlanarBdProp}, i.e., 
$\sum 	|B_{\rho}(\Cell_j)|\le K_0 \rho^2|H|^2$.
Then $K_1=2K_0$ satisfies the requirement of our Corollary, since otherwise
more than $(c/2)|H|^2$ cells with $\rho$--neighborhoods larger than
$(2K_0/c)\rho^2$ would give rise to a total strictly more than allowed.
\proofend
\medskip\par
%The second moment of such a subsystem can already
%be bounded by a quadratic function of $|H|$.
%\begin{cor}
%\label{SecondMomemtCor}
%From any  $c|H|^2$ cells one can select at least half of them
%for which 
%$$
%\sum 	|B_{\rho}(\Cell_j)|^2\le \frac{K_1^2}{c} \rho^4|H|^2.
%\proofend
%$$
%\end{cor}
\begin{rem}
The previous Corollary holds for $c|H|^3$ cells in $\R^3$, too, with an upper
bound 
$|B_{\rho}(\Cell_j)|\le({K_1}/{c})\rho^3 $ for another (larger) $K_1$.
\end{rem}
Unfortunately, neither Proposition~\ref{PlanarBdProp}, nor
Corollary~\ref{UniformBallsCor} is known in dimensions higher than three.
That is the main reason why we cannot extend Theorem~\ref{MainThm} to such
generality.
\end{subsection}

\begin{subsection}{The triangle selection lemma}
\begin{prop}
\label{EmoProp}
Let $\A(H)$ be a simple arrangement (i.e. the lines of $H$
are in general position).
Then for any $\Cell_j\in\A(H)$ and any $\rho\le|H|$ we have
$$
|B_{\rho}(\Cell_j)|>\frac{\rho^2}{32}.
$$
\end{prop}
\proof
We demonstrate the statement in two steps.
\medskip\pni
{\bf Step 1. }
$B_{\rho}(\Cell_j)$ contains at least ${\rho^2}/{8}$ vertices.\\
We copy the proof from \cite{Emo:arrangement}.
Draw a generic straight line which enters
$\Cell_j$ and select, on one of its two rays, the $\rho/2$ lines $h\in H$
which are closest to $\Cell_j$.
On each such line walk $\rho/2$ steps in any direction,
reaching at least $(\rho/2)\cdot(\rho/2)=\rho^2/4$ vertices
of $\A(H)$. Each of these is counted at most twice (since $\A$
is simple), whence the required inequality.
\medskip\pni
{\bf Step 2. }
Using the bound found in Step 1, we want to show
$|B_{\rho}(\Cell_j)|\ge{\rho^2}/{32}$.
To this end, consider the edges of $B_{\rho}(\Cell_j)$ as a 
graph on its vertex set. Denote by $v_4$ the number of vertices 
of degree exactly~4 and note that this set coincides 
with the vertex set of $B_{\rho-1}(\Cell_j)$.
Beyond them, there usually exist some other vertices along the
``outer boundary'', each of degree exactly~2.
If there is such an ``outer'' vertex, we delete it and 
glue together the two original edges incident on it to form one new
edge.
(This changes neither $v_4$ nor the number of faces.)
We repeat until only $v_4$ vertices remain, each
of degree~4. 
Then, using Euler's relation, 
$$
v_4+f=e+2=\frac{4v_4}{2}+2=2v_4+2,
$$
whence $f=v_4+2\ge (\rho-1)^2/8+2$,
by the lower bound in Step 1.
Thus
$$
|B_{\rho}(\Cell_j)|=\#\text{ of bounded faces } = f-1 
\ge \frac{(\rho-1)^2}{8}+1 >\frac{\rho^2}{32},
$$
since the inequality is obvious for $\rho<2$
and $\rho-1\ge\rho/2$ for $\rho\ge2$.\proofend
\begin{rem}
Proposition \ref{EmoProp} can be extended to $\R^d$ as 
$|B_{\rho}(\Cell_j)|>c_d{\rho^d}$, for some $c_d>0$,
independent of $\rho$.
\end{rem}
\begin{lem}[Triangle--selection Lemma]
\label{TriangSelLemma}
Let there be given a set $\P$ of 
$m$ points in $\R^2$ which lie in distinct cells
of a simple arrangement of $K m$ cells. 
($K$ is fixed while $m$ is large.)\\
If $\rho\ge\rho_0(K)=2048K+1024$,
then there exist at least m/6 triangles $U_1U_2U_3$, 
each spanned by three given points, such that
$\dist(U_iU_j)\le\rho$ for $1\le i,j\le3$.
\end{lem}
\proof
If a point $P_j$ in cell $\Cell_j$ is not the vertex of any such 
``small'' triangle then
$B_{\rho/2}(\Cell_j)\cap \P$ consists of collinear points
and thus 
$$
|B_{\rho/2}(\Cell_j)\cap\P|\le \rho+1 \le 2\rho.
$$
Consider the even smaller neighborhood $B_{\rho/4}(\Cell_j)$.
By Proposition~\ref{EmoProp}, it contains 
$$
\ge|B_{\rho/4}(\Cell_j)|-| B_{\rho/2}(\Cell_j)\cap\P|
    \ge\frac{\rho^2}{512}-2\rho
$$
empty cells.
Moreover, each such empty cell $\Cell_i$ was counted at most
$$
| B_{\rho/4}(\Cell_i)\cap\P|\le| B_{\rho/2}(\Cell_j)\cap\P|\le 2\rho
$$
times, giving a total of $\le K\cdot m\cdot2\rho$.
Thus the number of points which are incident upon no good triangle
is
$$
\le \frac{K\cdot m\cdot2\rho}{{\rho^2}/{512}-2\rho}=
   \frac{K\cdot m\cdot2}{{\rho}/{512}-2}=
   \frac{1024Km}{\rho-1024}\le\frac{m}{2},
$$
if $\rho\ge\rho_0(K)=2048K+1024$.
Therefore, at least $m/2$ points ARE incident upon at least one
good ``small'' triangle, giving
$$
\ge\frac{m/2}{3}=\frac{m}{6}
$$
distinct such triangles. \proofend
\begin{rem}\label{TriangSelRem}
In $\R^d$, if $\rho\ge\rho_0^{(d)}(K)$, then it is possible to select
$c_dm$ ``small'' simlpices (where $c_d>0$ and  $\rho_0^{(d)}(K)$
do not depend on $m$).
\end{rem}
\end{subsection}

\begin{subsection}{The ``Average--Forcing Lemma''}

We start with a simple observation.
\par
Let $\U_1$,  $\U_2$,  $\U_3$ be disjoint finite sets
and $\Delta$ a system of (ordered) triples $(U_1,U_2,U_3)$
such that $U_i\in\U_i$ for $i=1,2,3$.
As usual, we call the number of triples which contain
a given element $U_i$ the \emph{degree\/} of $U_i$ and
denote it by $\deg_{\Delta}(U_i)$.
We shall also use  the \emph{average degree in\/}  $\U_i$
which, of course, is $|\Delta|/|\U_i|$.
\begin{prop}[Folklore]
\label{FolkloreProp}
Given a triple system $\Delta$ as above,
there always exist subsets $\U_i\p\subset\U_i$ (for $i=1,2,3$)
such that the subsystem
$$
\Delta\p\defeq\{(U_1,U_2,U_3)\in\Delta\;;\;
    U_i\in\U_i\p   \text{\quad (for }i=1,2,3)\}
$$
satisfies
\begin{enumerate}[(i)]
\item
Each new degree is at least one quarter 
of the corresponding original average degree, i.e.,
for each $U_i\in\U_i\p$ we have
$$
\deg_{\Delta\p}U_i\ge\frac{1}{4}\cdot\frac{|\Delta|}{|\U_i|};
$$
\item
At least one quarter of the original triples are preserved, i.e.,
$$
|\Delta\p|\ge\frac{1}{4}\cdot\|\Delta|.
$$
\end{enumerate}
\end{prop}
\proof
Repeatedly delete those elements whose degree is less than required.
You cannot delete more than $3(|\Delta|/4)$ triples.
\proofend
\medskip\par
The setting of the following Lemma is ``2/3 Euclidean''
while ``1/3 abstract''.
We start with two (disjoint) point sets $\P_1,\P_2\subset\R^2$
and another (abstract) set $\S$ 
--- that is why the Lemma is  ``1/3 abstract''.
Moreover, let $\Delta$ be a triple system on these three sets,
as in the previous Proposition, with the additional requirement that
\par\centerline{\emph{
any pair $( S,P_i)$ is contained in at most one triple of $\Delta$.
			}}
\par\noindent
We assume that $|\P_i|\le N$ for $i=1,2$ while $|\S|=\lambda N$
and each $ S\in\S$ has $\deg_\Delta( S)\ge cN$,
for a fixed $c>0$.
Finally, let $H_1$ and $H_2$ be two systems of straight lines in $\R^2$
in general position, such that, for $i=1,2$,
all $P_i\in\P_i$ lie in distinct cells of $\A(H_i)$ while
 $|H_i|\le C\sqrt{N}$, for a fixed $C>0$.
\begin{lem}[Average--Forcing Lemma]
\label{AveForcingLem}
Beyond the setting outlined above,
let there be given an arbitrary $\rho>0$.
Then there exist
$\P_1^*\subset\P_1$, $\P_2^*\subset\P_2$,  and $\S^*\subset\S$ 
which, together with the corresponding subsystem $\Delta^*$, satisfy
\begin{enumerate}[(i)]
\item
$|\S^*|\ge c^*|\S|=c^*\lambda N$;
\item
for all $ S\in\S^*$ we have $\deg_{\Delta^*}S\ge c^*N$;
\item
consequently $|\P_i^*|\ge c^*N$,  for $i=1,2$;
\item
for any $P_i\in\P_i^*$ and the cell $\Cell_i\in\A(H_i)$ which
contains $P_i$, we have
$$
|B_{\rho}(\Cell_i)|\le C^*\rho^2;
$$
\end{enumerate}
where $c^*=c^*(c,C)$ and  $C^*=C^*(c,C)$ do not depend on $N$,
$\rho$, or $\lambda$.
\end{lem}
\proof
In  order to fulfill (iv) we shall, of course, use 
Corollary~\ref{UniformBallsCor} twice; once for $\P_1$ 
and once for $\P_2$.
However, at each step, we must make sure that the subsets selected
still participate in many triples of $\Delta$.
This will be achieved by always applying 
Proposition~\ref{FolkloreProp} before each such selection ---
and also at the very end, to guarantee (ii).
\par
Though this brief outline is certainly enough for the experienced reader,
we also describe all details as follows.
\medskip\pni
{\bf Step 1. }
First, using Proposition~\ref{FolkloreProp}, we select
$\P_1\p \subset\P_1$, $\P_2\p \subset\P_2$,   $\S\p \subset\S$ 
and the corresponding subsystem $\Delta\p $, such that
$|\Delta\p|\ge(1/4)|\Delta|\ge (1/4)c\lambda N^2$, whence
\[
\begin{aligned}
|\P_i\p|&\ge\frac{(1/4)c\lambda N^2}{\lambda N}= \frac{1}{4} c N =
\frac{c}{4C^2} \bigl(C\sqrt{N}\bigr)^2	\qquad\text{ and }\cr
|\S\p|&\ge\frac{|\Delta\p|}{|\P_i\p|}\ge\frac{(1/4)c\lambda N^2}{N}= 
\frac{1}{4} c\lambda N.\cr
\end{aligned}
\]
Moreover, for all $P_i\in\P_i\p$, 
\begin{equation}
\label{HighDegEq}
\deg_{\Delta\p}P_i \ge \frac{1}{4}\cdot\frac{c\lambda N^2}{N}= 
	\frac{1}{4}c\lambda N.
\end{equation}
%\pni
{\bf Step 2. } Using Corollary~\ref{UniformBallsCor} 
with $c/(4C^2)$ in place of $c$ there, we can select
a $\P_1\pp\subset\P_1\p$ of at least $(1/8)cN$ elements,
such that the cells ${\Cell_j\in\A(H_1)}$ which contain its points,
satisfy
$$
|B_{\rho}(\Cell_j)|\le K_1\frac{4C^2}{c} \rho^2.
$$
By (\ref{HighDegEq}), these $P_1\in\P_1\pp$ participate in a total of
$$
\ge\biggl(\frac{1}{8}c N\biggr)\cdot\biggl(\frac{1}{4}c\lambda N\biggr)=
 \frac{c^2}{32}\lambda N^2
$$
triples. We denote the set of these by $\Delta\pp$ while preserving
$\P_2\pp=\P_2\p$ and $\S\pp=\S\p$.
\par
At the moment we have
\[
\begin{aligned}
\frac{c}{\strut8}N\le& |\P_1\pp|\le N;\cr
\frac{c}{\strut4}N\le& |\P_2\pp|\le N;\cr
\frac{1}{4}c\lambda N\le &|\S\pp|\le \lambda N;\cr
&|\Delta\pp|\ge \frac{c^2}{32}\lambda N^2,\cr
\end{aligned}
\]
whence the average degree in $\S\pp$ is at least $({c^2}/{32}) N$.
\medskip\pni
{\bf Step 3. } As in Step 1, we use Proposition~\ref{FolkloreProp}
to select
$\P_1\ppp \subset\P_1\pp$, $\P_2\ppp \subset\P_2\pp$,   $\S\ppp \subset\S\pp$ 
and the corresponding subsystem $\Delta\ppp $, such that
$|\Delta\ppp|\ge(1/4)|\Delta\pp|\ge (c^2/128)\lambda N^2$, whence
\[
\begin{aligned}
|\P_i\ppp|&\ge\frac{(c^2/128)\lambda N^2}{\lambda N}=\frac{c^2}{\strut128}N=
	\frac{c^2}{128C^2} \bigl(C\sqrt{N}\bigr)^2\qquad\text{ and }\cr
|\S\ppp|&\ge\frac{|\Delta\ppp|}{|\P_i\ppp|}\ge\frac{(c^2/128)\lambda N^2}{N}=
	 \frac{c^2}{128}\lambda N.\cr
\end{aligned}
\]
Moreover, for all $P_i\in\P_i\ppp$, 
\begin{equation}
\label{SecondHighDegEq}
\deg_{\Delta\ppp}P_i \ge \frac{1}{4}\frac{(c^2/32)\lambda N^2}{N}= 
\frac{c^2}{128}\lambda N.
\end{equation}
%\pni
{\bf Step 4. } Using Corollary~\ref{UniformBallsCor} again, this time
with $c^2/(128C^2)$ in place of $c$ there, we can select
a $\P_2\pppp\subset\P_2\ppp$ of at least 
$(1/2)|\P_2\ppp|\ge (c^2/256)N$ elements,
such that the cells ${\Cell_j\in\A(H_2)}$ which contain its points,
satisfy
$$
|B_{\rho}(\Cell_j)|\le K_1\frac{128C^2}{c^2} \rho^2.
$$
By (\ref{SecondHighDegEq}), these $P_2\in\P_2\pppp$ participate in a total of
$$
\ge\biggl(\frac{c^2}{256} N\biggr)\cdot\biggl(\frac{c^2}{128}\lambda N\biggr)=
 \frac{c^4}{2^{15}}\lambda N^2
$$
triples. We denote the set of these by $\Delta\pppp$ while preserving
$\P_1\pppp=\P_1\ppp$ and $\S\pppp=\S\ppp$.
\par
At this moment we have
\[
\begin{aligned}
\frac{c^2}{\strut128}N\le& |\P_1\pppp|\le N;\cr
\frac{c^2}{\strut256}N\le& |\P_2\pppp|\le N;\cr
\frac{c^2}{128}\lambda N\le &|\S\pppp|\le \lambda N;\cr
&|\Delta\pppp|\ge \frac{c^4}{2^{15}}\lambda N^2,\cr
\end{aligned}
\]
whence the average degree in $\S\pppp$ is at least $({c^4}/{2^{15}}) N$
and at least $({c^4}/{2^{15}})\lambda N$ in the $\P_i\pppp$.
\pni
{\bf Step 5. } 
Finally, we use Proposition~\ref{FolkloreProp} for the third time to select
subsets 
$\P_1^*\subset\P_1\pppp$, $\P_2^*\subset\P_2\pppp$  and $\S^*\subset\S\pppp$ 
together with the corresponding subsystem $\Delta^*\subset\Delta\pppp$,
which already satisfy the requirements of Lemma~\ref{AveForcingLem}
with $c^*\defeq c^4/2^{16}$ and $C^*\defeq 128K_1C^2/c^2$. \proofend
\begin{rem}
\label{AveForcingRem}
Lemma~\ref{AveForcingLem} can be extended to $\R^3$ with 
$c^*\rho^3$ on the right hand side of (iv).
As for higher dimensions, no such result is known, due to the
lack of suitable generalizations of 
Proposition~\ref{PlanarBdProp} and Corollary~\ref{UniformBallsCor}.
\end{rem}
\end{subsection}

\begin{subsection}{Proof of the Main Theorem~\ref{MainThm}}

The graph of an affine mapping $\phi$ is a plane $S=S_{\phi}$ 
in $\R^2\times\R^2$.
The assumption of $|\phi(\P)\cap\P|\ge c|\P|$ is equivalent to
$|S_\phi\cap(\P\times\P)|\ge c|\P|$. In such a case we shall say that
the plane is ``$c|\P|$--rich''.
Thus it suffices to prove the following.
\begin{thm}
Let $\P_1,\P_2\subset\R^2$ with $|\P_1|=|\P_2|=N$ be proper 2--dimensional
point sets \utacfC.
Moreover, let $0<c<1$ be arbitrary.
Then
$$
\#\{\text{planes}\ S\subset\R^4 \;;\; 
      |S\cap(\P_1\times\P_2)|\ge cN\} \ \le\   C_1\cdot N,
$$
where $C_1=C_1(C,c)$ does not depend on $N=|\P_i|$.
\end{thm}
\proof
Denote the set of $cN$--rich planes by $\S$ and let
$|\S|=\lambda N$, for some positive $\lambda$.
We want to bound  $\lambda$ from above by a constant $C_1$.
To this end we shall use Lemmata~\ref{TriangSelLemma} and
\ref{AveForcingLem}  for a value of $\rho=\rho(C,c)$,
to be specified later  in terms of
$c^*=c^*(c,C)$ and  $C^*=C^*(c,C)$, given in 
Lemma~\ref{AveForcingLem}.
\medskip\par
We shall call the planes $\{(x,y,0,0)\;;\;x,y\in\R \}$
and $\{(0,0,z,w)\;;\;z,w\in\R \}$ in $\R^4$ 
--- which contain the copies of the $\P_i$ generating $\P_1\times\P_2$
--- the first and second coordinate plane, respectively.
\medskip\par
First we define the triple system
$$
\Delta\defeq
       \bigl\{(P_1,P_2,S)\;;\;P_i\in\P_i,S\in\S,\{P_1\}\times\{P_2\}\in S  
       \bigr\}
$$
and apply Lemma~\ref{AveForcingLem} to $\Delta$ and the two
systems $H_1$, $H_2$ of straight lines which cut the $\P_i$
into singletons, according to the ``proper 2D'' assumption.
Thus we can find subsets
$\P_1^*\subset\P_1$, $\P_2^*\subset\P_2$,  and $\S^*\subset\S$ 
together with the corresponding subsystem $\Delta^*$,
with the following properties.
\begin{enumerate}[(a)]
\item
$|\S^*|\ge c^*|\S|=c^*\lambda N$ and
for all $ S\in\S^*$ we have $\deg_{\Delta^*}(S)\ge c^*N$.
For these $S$, we introduce 
\begin{enumerate}[(i)]
\item
the notation 
$m_\S\defeq|S\cap(\P_1^*\times\P_2^*)|\ge c^*N$;
\item
and also a new arrangement $\A_S$ defined by the union of
projections (to $S$) of $H_1$ and $H_2$ (considered 
in the first and second coordinate plane and
cutting $\P_1$ and $\P_2$, respectively, into singletons).
\end{enumerate}
Thus  
$$
|\A_S|\le(2C\sqrt{N})^2=4C^2N=\frac{4C^2}{c^*}\cdot c^*N \le
 \frac{4C^2}{c^*}\cdot m_S.
$$
Putting $K={4C^2}/{c^*}$ and defining $\rho=\rho(K)=2048K+1024$
we can apply Lemma~\ref{TriangSelLemma}, to find  in each $S\in\S^*$
at least $m/6\ge c^*N/6$ distinct ``small'' triangles whose vertices
are in $\P_1^*\times\P_2^*$ and any pair of vertices 
has $\dist\le\rho$ as measured in $\A_S$.
\item
On the other hand, for $i=1,2$, 
each $P_i\in\P_i^*$ lies in an original cell 
$\Cell_i\in\A(H_i)$ with $|B_{\rho}(\Cell_i)|\le C^*\rho^2$.
Therefore, the number of small triangles which have $P_i$
as a vertex, cannot exceed $|B_{\rho}(\Cell_i)|^2\le (C^*)^2\rho^4$.
(Note that the distance measured in $\A_S$ 
is at least as much as if measured in any $\A(H_i)$.)
Hence the total number of such triangles is 
$\le (C^*)^2\rho^4|\P_1^*||\P_2^*|=(C^*)^2\rho^4N^2$.
\end{enumerate}
Putting (a) and (b) together for the already defined $\rho$, we have
$$
(c^*\lambda N)\cdot(c^*N/6)\le \#\text{ of small triangles } \le
  (C^*)^2\rho^4N^2,
$$
whence
$$
\lambda\le \frac{6(C^*)^2}{(c^*)^2}\cdot\rho^4\defeq C_1(C,c),
$$
a constant in terms of the parameters $C$ and $c$.\proofend
\begin{rem}
\label{ThreeDimProofRem}
Essentially the same idea works in $\R^3$, with Remarks~\ref{TriangSelRem}
and \ref{AveForcingRem} in place of Lemmata~\ref{TriangSelLemma}
and \ref{AveForcingLem}, respectively.
\end{rem}
\begin{rem}
An affirmative answer to Problem~\ref{RefinementProbl} could substitute 
the very technical and
complicated Lemma~\ref{AveForcingLem} in the proof of Theorem~\ref{MainThm}
as follows.
We first refine $H_1$ and $H_2$, by adding $O(\sqrt{N})$ lines,
to get some $H_1^+$ and $H_2^+$ with small $\rho$--neighborhoods.
Then we use the Triangle--Selection Lemma~\ref{TriangSelLemma} for the
corresponding $\A_S^+$ and finish as in Step (b).
\end{rem}
This could even extend to higher dimensions, provided that
Problem~\ref{RefinementProbl} can.
\end{subsection}
\end{section}

\bibliographystyle{alpha}%{plain}

\begin{thebibliography}{Wel92}

\bibitem[EK01]{EGyKZ}
Gy{\"o}rgy Elekes and Zolt{\'a}n Kir{\'a}ly.
\newblock On combinatorics of projective mappings.
\newblock {\em Journal of Algebraic Combinatorics}, 14:183--197, 2001.

\bibitem[Ele97]{EGy:LinI}
Gy{\"o}rgy Elekes.
\newblock On linear combinatorics {I}.
\newblock {\em Combinatorica}, 17 (4):447--458, 1997.

\bibitem[Ele98]{EGY:LinII}
Gy{\"o}rgy Elekes.
\newblock On linear combinatorics {II}.
\newblock {\em Combinatorica}, 18 (1):13--25, 1998.

\bibitem[Ele02]{EGy:survey}
Gy\"orgy Elekes.
\newblock Sums versus products in number theory, algebra and {Erd\H{o}s}
  geometry --- a survey.
\newblock {\em in: {\sl Paul Erd\H{o}s and his Mathematics II} Bolyai Math.
  Soc. Stud {\bf 11}}, pages 241--290, 2002.

\bibitem[Ele05]{EGy:DoFI}
Gy\"orgy Elekes.
\newblock On the dimension of finite point sets {I}.
\newblock el\H{o}k\'esz\"uletben, 2005.

\bibitem[ER00]{EGyRL}
Gy{\"o}rgy Elekes and Lajos R{\'o}nyai.
\newblock A combinatorial problem on polynomials and rational functions.
\newblock {\em Journal of Combinatorial Theory, series A}, 89:1--20, 2000.

\bibitem[PS98]{PaSh:98}
J\'anos Pach and Micha Sharir.
\newblock On the number of incidences between points and curves.
\newblock {\em Combinatorics, Probability and Computing}, 7:121--127, 1998.

\bibitem[ST83]{SzTr:83}
Endre Szemer\'edi and W.~T. T{rotter~Jr.}
\newblock Extremal problems in {D}iscrete {G}eometry.
\newblock {\em Combinatorica}, 3 {\rm (3--4)}:381--392, 1983.

\bibitem[Wel92]{Emo:arrangement}
Emo Welzl.
\newblock On spanning trees with low crossing numbers.
\newblock {\em ``Data Structures and Efficient Algorithms'', (B. Monien, Th.
  Ottmann, Eds.), Lecture Notes in Computer Science}, 594:233--249, 1992.

\end{thebibliography}

\end{document}